\documentclass[12pt]{article}

\usepackage{amssymb, latexsym, amsmath}
\newtheorem{Theorem}{Theorem}[section]
\newtheorem{Proposition}[Theorem]{Proposition}
\newtheorem{Lemma}[Theorem]{Lemma}
\newtheorem{Corollary}[Theorem]{Corollary}
\newtheorem{Remark}[Theorem]{Remark}

\newcommand{\bi}{\begin{enumerate}}
\newcommand{\ei}{\end{enumerate}}
\newcommand{\be}{\begin{equation}}
\newcommand{\ee}{\end{equation}}
\newcommand{\ba}{\begin{array}}
\newcommand{\ea}{\end{array}}

\def\c{\mathbb{C}}
\def\cp{\mathbb{CP}}
\def\r{\mathbb{R}}

\def\z{\mathbb{Z}}

\def\lck{l.c.K. }
\def\T{\mathbf{T}}
\def\t{\mathbf{t}}

\author{A.Fujiki, M.Pontecorvo}
\title{Bi-Hermitian metrics on Kato surfaces }
\begin{document}
\maketitle 

\begin{abstract}
We produce the first examples of bi-Hermitian metrics with connected numerically 
anti-canonical divisor on Kato surfaces. 
The picture of compact complex surfaces admitting such structures remains open 
only for intermediate Kato surfaces.
\end{abstract}

\section{Introduction}
By a bi-Hermitian surface we will mean in this work a Riemannian four--manifold $(M^4,g)$
with two integrable complex structures $J_+$ and $J_-$ which are $g$-orthogonal
and induce the same orientation. 
To avoid the trivial case we assume there is $x\in M$ where $J_+(x) \neq \pm J_-(x)$.
The compatibility of $J_{\pm}$ with a Riemannian metric -- hence conformal structure $c:=[g]$ --
and an orientation is equivalent to the algebraic condition that, viewed as sections of $End(TM)$,
they satisfy the following algebraic relation for some smooth function $p:M\to [-1,1]$,  
see \cite{tr74} and \cite[1.5]{po97}         $$J_+J_- + J_-J_+ = -2p\cdot id.$$

The four-dimensional manifold $M$ will always be assumed to be compact and in this
situation it was shown in \cite[1.6]{po97} that a bi-Hermitian surface $(M,c,J_\pm)$ can admit one more compatible complex
structure if and only if it is actually hyperhermitian: therefore a $K3$ surface, a complex torus or a
special diagonal Hopf surface, all with standard metrics, by Boyer classification  \cite{bo88}.

Sort of conversely, one can consider the problem of finding all compact complex surfaces 
$S=(M,J)$ which admit a Hermitian metric $g$ compatible with a different complex structure on $S$.
The underlying conformal four--manifold $(M,c,J_{\pm})$ is bi-Hermitian and
the complex surface $S$ is either $(M,J_+)$ or $(M,J_-)$. 
Examples show that these two surfaces may or may not be biholomorphic or deformation equivalent.

The search for  bi-Hermitian surfaces is of particular interest because in higher dimension there is more flexibility
for compatible complex structures and the above mentioned result \cite[1.6]{po97} certainly does not hold \cite{sa09}.

\medskip
The first obstruction to existence of bi-Hermitian structures which was observed \cite{po97}, \cite[Prop.2]{agg99} is that such a surface $S$ needs to have an (automatically effective, or zero) numerical anticanonical divisor
\be    \T=F\otimes K^*   \label{fund}  \ee
where $K$ will always denote the canonical line bundle,
$K^*=K^{-1}$ the anticanonical bundle and $F$ is a holomorphic line bundle of 
zero Chern class which is always either trivial or real of negative degree. 
Using the terminology of Dloussky \cite{dl06} we will refer to $\T$ as a NAC 
(numerically anti-canonical) divisor.

In this work we will always identify, by abuse, a divisor with the corresponding holomorphic
line bundle. Although in general $\T$ is highly non-reduced,
set theoretically its support is the set where the smooth function $p$ has absolute value $1$;
in other words the support of $\T$ is the set where the two complex structures agree up sign:
$$ \operatorname{supp} \T=\{x\in M : J_+(x)=J_-(x)\} \amalg \{x\in M : J_+(x)=-J_-(x)\} .$$
 
\begin{Remark} \rm It is worthwhile noticing that each of the above two components turns out to be either empty or connected - by surface classification - so that if we set $\mathbf t :=b_0(\T)$ there are only three possibilities:  $\mathbf t =0,1,2$ all of which do actually occur.
\hfill $\bigtriangleup$ 
\end{Remark}

\smallskip The case $b_1(S)$ even is well understood: the line bundle $F$ of equation (\ref{fund}) turns out to be always trivial: $F=\mathcal O _S$ by \cite[Lemma 4]{agg99}. After preliminary results of Hitchin on rational surfaces \cite{hi06} and  \cite{hi07}, Goto proved in \cite{go12} that a compact complex surface $S$ with even first Betti number admits bi-Hermitian metrics if and only if $H^0(S,K^*)\neq 0$  showing that (\ref{fund}) is in fact a sufficient condition for surfaces of K\"ahler type.

We thus concentrate on the case $b_1(S)$ odd where it was shown by 
Apostolov \cite{ap01} that the fundamental equation (\ref{fund})
forces such a surface to be in class VII so that $b_1(S)=1$ always. 
Furthermore, since blowing up points is compatible with bi-Hermitian metrics \cite{fp10} and \cite{cg11}, we will also assume that $S$ is minimal -- i.e. free from smooth $\mathbb{CP}_1$'s of self-intersection $-1$. Under this minimality assumption it was then shown that $S$ must be a Hopf surface  when $b_2(S)=0$ \cite{ap01} or a  Kato surface when $b_2(S)>0$ \cite{dl06}.

\medskip
We present in this work a new result about bi-Hermitian Kato surfaces with $\t=1$, see 
Theorem \ref{parinoue}  which is complementary to the construction 
in \cite{fp10} and leads to a fairly complete picture for  ``unbranched" Kato surfaces.
As the Hopf surface case was well understood in \cite{po97}, \cite{ad08} and \cite{abd15}
we obtain a classification which shows a surprising abundance of bi-Hermitian surfaces,
except for the case of ``intermediate" Kato surfaces which, contrary to what was announced in an early version of this work \cite{fp16}, remains open.

\medskip
Our main source of inspiration is the work of \cite{abd15} which points out how the $b_1(S)$-odd
case should be thought of as a twisted K\"ahlerian case, the twisting being given by the flat bundle $F$
in (\ref{fund}). For this reason we now begin with a brief discussion of Kato surfaces and special flat line bundles. 

\medskip  \noindent 

\section{Kato surfaces and flat holomorphic line bundles of real type}

As we mentioned in the introduction, when looking for compact complex surfaces $S$ admitting a bi-Herimtian structure the only open case is that of Kato surfaces. They are compact complex $2$-manifolds belonging to class VII$_0$ in Kodaira's classification, meaning that the first Betti number $b_1(S)$  is odd, the Kodaira dimension Kod$(S)=-\infty$ and that every smooth rational curve $C\subset S$ has self-intersection $C^2\neq -1$ in fact $C^2\leq -2$, otherwise $S$ would be a rational surface. Surfaces in this class always have first Betti number $b_1(S)=1$ and they are classified only when the second Betti number $b_2(S)=0$ in which case $S$ is either a Hopf surface -- i.e. the universal covering is $\c^2\setminus\{0\}$ -- or else is a Inoue-Bombieri surface; while bi-Hermitian Hopf surfacs are well understood from the work of Apostolov-Dloussky \cite{ad08}, no bi-Hermitian metrics can be supported on Inoue-Bombieri surfaces \cite{ap01}.

When $b_2(S)>0$ we say that $S\in$VII$_0^+$; the only known examples of such surfaces are called Kato surfaces \cite{ka77} and contain a \it global spherical shell \rm (GSS) meaning that there is a holomorphic embedding $\phi:U\to S$, where $U\subset\c^2\setminus\{0\}$ is a neighborhood of the unit sphere $S^3\subset\c^2$, such that $S\setminus \phi(U)$ is connected.

Such a surface $S$ can always be obtained by successive blow ups of the unit ball in $\c^2$ starting from the origin and continuing to blow up a point on the last constructed exceptional curve until, in order to get a compact complex surface, we take the quotient by a biholomorphism $\sigma$ which identifies a small ball centered at a point $p$ in the last exceptional curve with a small ball around the origin $0\in\c^2$ in such a way that $\sigma(p)=0$, see the picture on \cite[p.5]{dl14}. This process constructs $b:=b_2(S)$ rational curves $D_1,\dots , D_b$ of self-intersection $\leq -2$ which are the only rational curves in $S$ and is therefore minimal. A strong geometric characterization of these surfaces has been given by Dloussky-Oelejklaus-Toma who showed that $S\in$VII$_0^+$ is a Kato surface if and only if it admits $b=b_2(S)$ rational curves \cite{dot03}. It should also be mentioned that the global spherical shell conjecture \cite{na89} states that every $S\in$VII$_0^+$ has a GSS and is therefore a Kato surface.  We would now like to give a quick description of this class of non-K\"ahler surfaces based on the configuration of their rational curves using the terminology of the Japanese school.

When every blow up occurs at a generic point of the previously created exceptional curve we get $D_i^2=-2$ for each $i=1,\dots , b$ and the divisor $D_1+ \cdots + D_b$ is a cycle $C$ of rational curves of self-intersection $C^2=0$. In this particular case $S$ is called an \it Enoki \rm surface and is the only case in which $D_1, \dots , D_b$ do not span $H^2(S,\z)$.
An Enoki surface admitting a compact curve $E$ which is not rational smooth is called a \it parabolic Inoue \rm surface, in this case $E$ must be smooth elliptic with $E^2=-b$; these surfaces always admit a holomorphic vector field and are the only Enoki surfaces with these properties.

The other extreme case is that of each blow up point being at the intersection with a previously created exceptional curve; these surfaces are called \rm hyperbolic Inoue \rm when the maximal curve is the union of two cycles and the only other case is that of their $\z_2$-quotients called \it half-Inoue \rm surfaces. All other surfaces are obtained by a mixed procedure in which there are generic blow ups as well as blow ups at intersection points and are therefore called \it intermediate \rm Kato surfaces.

Contrary to an early version of this work \cite{fp16} we don't know about existence of bi-Hermitian metrics in the intermediate case. The situation for the  ``extreme" Kato surfaces is more clear because anti-self-dual bi-Hermitian metrics were constructued in \cite{fp10} on hyperbolic and parabolic Inoue surfaces, while half Inoue and the other Enoki surfaces cannot admit such structures at all \cite{fp14} \cite{abd15}.
The present work concludes the classification of ``extreme" Kato surfaces admitting bi-Hermitian metrics by producing the only remaining case of parabolic Inoue surfaces with $\t=1$. 
{\bf Notation.} The dual graph of all rational curves on intermediate Kato surfaces always consists of one cycle with a positive number of trees attached to it, see the picture in \cite{fp15} and we will also call them ``branched" Kato surfaces for this reason. In the extreme case there are only cycles with no trees and we will call hem ``unbranched".

\medskip  \noindent 
\it Flat holomorphic line bundles of real type  \rm

Denoting by $(S,J)$ or simply by $S$ either of the two compact complex surfaces $(M,J_{\pm})$
we will be particularly interested in $\operatorname{Pic}^0_\r(S)$ by which we mean holomorphic 
line bundles in the image of $H^1_{dR}(S,\r)$ and in the kernel of the first Chern class map.
Recall that for class-VII surfaces we have an isomorphism 
$H^1(S,\mathcal O)\cong H^1(S,\c)\cong \c$, the group of zero-Chern class holomorphic line bundles
is non-compact $\operatorname{Pic}^0(S)\cong H^1(S,\c^*)\cong \c^*$ and we will denote by 
$\operatorname{Pic}^0_\r(S)$ the bundles of real type coming from the exponential map
$H^1_{dR}(S,\r)  \xrightarrow{\exp} H^1(M,\r^{>0}) \hookrightarrow H^1(S,\c^*)$  
so that elements of $\operatorname{Pic}^0_\r(S)$ are complexifications of line bundles whose transition functions are real and locally constant.

In details, a closed $1$-form $\alpha$ is locally exact: $\alpha=df_i$ on the open cover $\{U_i\}$ it defines a line bundle $L_\alpha$ with positive local trivializations $e ^{f_i}$ on $U_i$ and locally constant transition functions $e^{f_{j}-f_{i}}$ on $U_i\cup U_j$ so that $L_\alpha\in H^1(S,\r^{>0})$ has vanishing first Chern class. 

This correspondence is compatible with cohomolgy and provides the following sequence of isomorphisms:
$$ \exp : H^1(S,\r) \to H^1(S,\r^{>0}) \to \operatorname{Pic}^0_\r(S) \hookrightarrow \operatorname{Pic}^0(S)\cong \c^*$$
We will use Convention 2.1 in \cite{ad16} and identify two real flat line bundles $L_\alpha$
and $L_{\alpha^{'}}$ in $H^1(S,\r^{>0})$ if they correspond to closed $1$-forms in the same 
cohomology class $a\in H^1_{dR}(S,\r)$.
Furthermore, we will denote by $\mathcal L  _a \in  \operatorname{Pic}^0_\r(S)$ 
the flat holomorphic line bundle of real type which is the complexifiation of 
$L_a \in H^1(S,\r^{>0}) = \exp( H^1(S,\r))$.

\medskip
We can now describe two main reasons for our interest in flat holomorphic line
bundles of real type and we will denote the correspondence above by the morphism
$$ \ba    {cccc}
 \exp: & H^1_{dR}(S,\r) &  \longrightarrow &  \operatorname{Pic}^0_\r(S)\cong \r^{>0} \\
                 & a                        &  \longmapsto &  \mathcal L _a
\ea  \label{realflat} $$
which sends $a+b$ to the tensor product 
$\mathcal L _a \otimes \mathcal L _b = \mathcal L _{(a+b)}$. 
While for the dual line bundle we use the notation $\mathcal L _a^* = \mathcal L _{-a}$.

We start by introducing the \it Lee form \rm of an arbitrary Hermitian metric $g$
on the complex surface $S$. In real dimension four the fundamental $(1,1)$-form
$\omega(\cdot , \cdot ) := g(J\cdot,\cdot)$ being non-degenerate induces an isomorphism
$\omega: \Lambda^1 \to \Lambda^3$ taking  $a\mapsto a\wedge\omega$ 
therefore there always exists $\theta$ such that $d\omega = \theta\wedge\omega$.
This real $1$-form $\theta$ is called the Lee form of $(g,J)$ and notice that $b_1(S)=odd$ 
implies that its harmonic part $\theta_h\neq 0$ by a result of Gauduchon, see \cite[Prop.1]{ag07}. 
Under conformal rescaling of the metric $e^h g$ the Lee form transforms as $\theta + dh$, so that its de Rham class remains well defined.

In order to obtain flat line bundles we need to look at \it closed \rm $1$-forms.
The first source comes from locally conformal K\"ahler metrics which
we abbreviate by l.c.K. and correspond to the condition $d\theta=0$ -- 
i.e. $\theta=df$ locally so that $e^{-f} g$ is a local K\"ahler metric on $S$ conformal to the original one.
We get in this way a real flat bundle $L_\theta$ whose complexification has the structure
of a  holomorphic flat line bundle of real type $\mathcal L  _\theta$ with operator
$\overline{\partial}_\theta := \overline{\partial} - \theta^{0,1}\wedge$.
We call this bundle the \it Lee bundle \rm of the l.c.K. metric; it is never trivial because $b_1(S)$ is odd.
It also follows from the above discussion that the \lck fundamental $(1,1)$-form $\omega$ 
is a $2$-form with values in $L_\theta ^*$, 
closed with respect to the differential operator $d_\theta:=(d- \theta\wedge)$.

A second source of closed real $1$-forms on a compact complex surface
comes from bi-Hermitian structures  $(M,c,J_\pm)$. For any $g\in c$ consider
the two Lee forms $\theta_\pm$ of the Hermitian metrics $(g,J_\pm)$ which of
course are real and need not be closed individually; 
however by \cite[Lemma 1]{agg99} the sum $\theta_+ + \theta_-$ is always closed and   
it was shown in \cite[Prop. 2]{agg99} that there always is a holomorphic section of
$F\otimes K^*$ where $F$ is the flat holomorphic line bundle of real type 
$F:=\mathcal L_{\frac12 (\theta_+ + \theta_-)}$.

The resulting effective (or zero) divisor $\T:=F\otimes K^*$ is then a NAC divisor in both 
$(M,J_{\pm})$ and recall that $0\leq \mathbf t :=b_0(\T) \leq 2$ by surface classification.

The case $\mathbf t =0$ was completely described in \cite{ad08} while it is shown in \cite{fp14} 
that $\mathbf t =2$ corresponds to the so called  \it generalized K\"ahler \rm case
-- i.e.  $\theta_+ = - \theta_ -$ \cite{ag07}.
In this work we will mainly be concerned with the situation when $\mathbf t =1$ and
$S$ is a Kato surface -- i.e. a minimal complex surface with a global spherical shell
and $b_2(S)>0$.

\bigskip\bigskip  \noindent
\it Degree of holomorphic line bundles \rm

By a strong result of Gauduchon \cite{ga84}, in every conformal class $c$ of Hermitian metrics
on a compact complex surface $S$ there is a unique representative $g\in c$ whose
Lee form $\theta$ is coclosed or equivalently whose fundamental $(1,1)$-form $\omega$
is $\partial\overline\partial$-closed.  
Such a $g$ will be called a  \it Gauduchon metric; \rm it can be used to define the notion of degree 
of a holomorphic line bundle $\mathcal L \in \operatorname{Pic}(S)$ as follows \cite{ga81}:
$$ \operatorname{deg}_g(\mathcal L) = \frac{1}{2\pi} \int_M \rho\wedge\omega $$
where $\rho$ is the real curvature of a holomorphic connection on $\mathcal L$,
it is a closed $2$-form  well defined modulo $\partial\overline\partial$-exact forms.
In this way the degree measures the $g$-volume with sign of a virtual meromorphic section
of $\mathcal L$: in particular $deg (\mathcal L)>0$ whenever 
$H^0(S,\mathcal L)\neq 0$ and $\mathcal L \neq \mathcal O$.

It is shown in \cite{ad16} that in the special case $\mathcal L$ is a flat holomorphic line bundle 
of real type so that $\mathcal L = \mathcal L _a$ for some $a\in H^1_{dR}(M,\r)$,  then: 
   \be 2\pi \operatorname{deg}(\mathcal L_a)=-\langle a,\theta \rangle = 
        -\langle a_h,\theta_h \rangle 
     \label{scalarproduct}  \ee
where $\theta$ is the Lee form of the Gauduchon metric $g$ with $g$-harmonic part denoted by
$\theta_h$ while $a_h$ denotes the $g$-harmonic part of the $1$-form $a$ and
$\langle , \rangle$ is the global inner product with respect to $g$; 
the last equality holds because $a$ is closed and $\theta$ is coclosed.

In our case we always have $b_1(S)=1$ and therefore the composition
\be   \operatorname{deg}_g : H^1_{dR}(S,\r) \xrightarrow{\exp} \operatorname{Pic}^0_\r(S) \to\r \ee
is a linear isomorphism since $\theta_h\neq 0$.

The isomorphism itself depends on the Gauduchon metric, however because the space of
Gauduchon metrics is connected (actually convex) and $H^1_{dR}(S,\r)\setminus \{0\}$ is not,
it follows that --  for a fixed de Rham class $a$ --
the sign of the degree of $\mathcal L _a \in \operatorname{Pic}^0_\r(S)$
does not depend on the particular conformal class but depends on the complex structure $J$ alone.

Following \cite{abd15},  in the case $b_1(S)=1$ one can therefore define an orientation on
$H^1_{dR}(S,\r)\cong\r$ by setting $a>b$ if and only if 
$\operatorname{deg} (\mathcal L _a) > \operatorname{deg} (\mathcal L _b) $ 
for some (hence any) Gauduchon metric on $(S,J)$.
 
In the present b-Hermitian context $(M,J_\pm)$, for a fixed class $a\in H^1_{dR}(M,\r)$ the sign of the degree of $\mathcal L _a$
may actually depend on the particular complex structure $J_+$ or $J_-$  as the following proposition shows.

Before getting to this result we recall a link with generalized K\"ahler geometry of Gualtieri 
who showed that this geometry always gives rise to bi-Hermitian metrics satisfying certain
equations which in our four-dimensional case $(M,c,J_\pm)$
amount to the condition that the two Lee forms admit a common  Gauduchon metric in $c$
and satisfy $\theta_+ =-\theta_-$ with respect to this metric. 
To be more precise there actually is a twisting coefficient $[H]\in H^3(M,\r)$ which turns out to vanish if and only if $b_1(M)$ is even \cite[Prop. 3]{ag07}. 

Therefore in the case under study -- i.e. $b_1(M)=$ odd -- 
a \it twisted generalized K\"ahler \rm structure is a bi-Hermitan surface with trivial fundamental flat line bundle: $F=\mathcal O$ because $F= \mathcal L_{\frac12 (\theta_+ + \theta_-)}$ by \cite{agg99}.

We can now state our result as follows; it will be used in the proof of \ref{real} and \ref{blowup}.

\begin{Proposition}  \label{minus}
A compact bi-Hermitian surface $(M,c,J_\pm)$ with $b_1(M)=1$ is twisted generalized K\"ahler
if and only if $\mathbf t =2$, if and only if the two complex structures $J_+$ and $J_-$
 induce opposite orientations on $H^1_{dR}(M,\r)$.
\end{Proposition}
 {\it Proof.} The fact that $\mathbf t =2$ is equivalent to $F=\mathcal O$ was shown in one direction in
 \cite[Prop. 4]{agg99} and in \cite[4.7]{fp14} for the other direction.
 
 Now, the degree of $F$ was computed in \cite{ap01} using the already mentioned 
 result of \cite{agg99} that the fundamental line bundle satisfies 
 $F=\mathcal L _{\frac12 (\theta_+ +\theta_-)}$.
It then follows from equation (\ref{scalarproduct}) that the degree with respect to the 
$J_+$--Gauduchon metric is
 $2\pi \operatorname{deg}_{J_{+}} (F) = -\langle \theta_+ , \frac12 (\theta_+ +\theta_-) \rangle =
 -\frac14 ||\theta_+ +\theta_- ||^2$  because the equality $||\theta_+||^2 = ||\theta_- ||^2$ 
 holds for any metric in the conformal class \cite{po97}.
 
 The first consequence is that $F=\mathcal O$ if and only if $\t=2$, if and only if 
  $\theta_+ = - \theta_- $.
But in this case   $\delta\theta_- =0$ with respect to the $J_+$-Gauduchon metric
and by equation (\ref{scalarproduct}) we have that for any $a\in H^1_{dR}(,\r)$
$$ 2\pi \operatorname{deg}_{J_{+}} (\mathcal L _a)  = -\langle \theta_+ , a \rangle =
   \langle \theta_- , a \rangle = - 2\pi \operatorname{deg}_{J_{-}}(\mathcal L _a) $$
  with respect to the common $J_\pm$-Gauduchon metric, showing that the induced
  orientations are opposite.
  
Finally, when $\t \neq 2$ we compute 
$2\pi \operatorname{deg}_{J_{-}}(F)=-\langle \theta_- , \frac12 (\theta_+ +\theta_-) \rangle
=-\frac14 ||\theta_+ +\theta_- ||^2 < 0$  with respect to any $J_-$ -Gauduchon metric,
as $F\neq \mathcal O$.
We conclude that the fundamental line bundle $F$ has negative degree for both $J_+$ and $J_-$ and this is enough to infer that the induced orientations agree in this case, as $b_1 (M)=1$.
 \hfill$\Box$\medskip

We end this section by observing that the degree of a Lee bundle -- i.e. the holomorphic flat line bundle of a l.c.K. metric -- is always strictly negative, because from equation (\ref{scalarproduct}) 
it follows that $\operatorname{deg} (\mathcal L _{\theta}) = -||\theta_h||^2 < 0$ otherwise $g$ 
would be conformally K\"ahler and $b_1(S)$ even. 
In particular then, $H^0(S,\mathcal L _{\theta})=0$ and the same holds for any positive power.

\section{Anti-self-dual bi-Hermitian surfaces}

A twistor approach to compact bi-Hermitian surfaces appeared in \cite{po97} and was
further developed in \cite{fp10}.
The main results are the following:  assume $(M,c,J_\pm)$ is compact and bi-Hermitian
with $b_1(M)= odd$; then the conformal class of the metric is anti-self-dual if and only if it is \lck with respect to a common Gauduchon metric. 
In this situation the fundamental divisor $\T=K^*$ must be anticanonical with two connected 
components -- i.e. $\t=2$.
The two Lee forms satisfy the generalized K\"ahler equation $\theta_+ = -\theta_-$.
In particular $\T$ is effective and both surfaces $S_\pm$ have negative Kodaira dimension.
Viceversa, $\t=2$ implies that $b_1(M)$ is odd \cite[Prop. 4]{agg99} as well as the generalized
K\"ahler condition \cite[Prop. 4.7]{fp14}.

The fact that the anticanonical bundle is effective and disconnected
has strong implications on the surfaces $S_\pm$: by \cite[Lemma 3.2]{fp10}
the complex structure -- recall once again that we assume the surface is minimal --
can only be that of a Hopf surface, a parabolic Inoue surface
or a hyperbolic Inoue surface. And the main contribution of \cite{fp10} is to show that
conversely all these surfaces do admit bi-Hermitian metrics, provided
a certain ``reality" condition is satisfied.

The aim of this section is to compute the Lee bundle of anti-self-dual bi-Hermitian surfaces,
clarify the reality condition in all of the possible three cases and derive some consequences.

\medskip
Before stating our result we need to set up the notation and recall some facts.
We will write $C+E$ for the anticanonical divisor $K^*$ of an anti-self-dual bi-Hermitian surface
$S$ where $C$ and $E$ are its connected components. $C$ will be a cycle of rational curve
in the case of parabolic or hyperbolic Inoue surface and a smooth elliptic curve in the remaining
case of Hopf surfaces. While $E$ will be smooth elliptic for Hopf and parabolic Inoue surfaces
and a cycle of rational curves for hyperbolic Inoue surfaces \cite{fp10}.

As the conformal class of the metric is anti-self-dual the total space of the twistor fibration
$t:Z\to M$ is a complex three-manifold \cite{ahs},
each compatible complex structure $J$ defines a smooth and disconnected divisor
$X=\Sigma + \bar{\Sigma}\subset Z$ which is the image of the tautological sections
$J$ and $-J$ of the twistor fibration. In the bi-Hermitian case we get four such sections
denoted by $X_\pm = \Sigma_\pm +  \bar{\Sigma}_\pm$; the sum
$X_+ + X_-$ represents an anticanonical divisor in the twistor space $Z$
because the Lee forms satisfy $\theta_+ + \theta_- =0$  \cite{po97}.

We also recall from \cite{po97} and \cite{fp10} that the two twistor divisors
$X_+$ and $X_-$ intersect each other --
transversally, except that they are tangential at the nodes of $C$ and $E$ --
in the following configuration:
the intersection $X_+ \cap X_-$ has four connected components consisting of
the maximal curve $C_+ +E_+=\Sigma_+\cap X_-$ in $\Sigma_+ \cong (M,J_+)$ and
the maximal curve $C_- +E_-=\Sigma_-\cap X_+$ in  $\Sigma_- \cong (M,J_-)$;
both maximal curves turn out to be anticanonical and reduced with two connected components..

We can arrange labelings in such a way that the cycle of rational curves is the complete intersection $C_+=\Sigma_+ \cap \Sigma_-$. 
By restricting to $C_+\subset\Sigma_+$ a locally defining function for the smooth divisor $\Sigma_-$ in $Z$ we get that 
\be
\nu_{C_{+} / \Sigma_{+}}=\Sigma_{-_{ |_{C_{+}}}} \quad \text{ and similarly } \quad
   \nu_{C_{+} / \Sigma_{-}}=\Sigma_{+_{ |_{C_{+}}}}   \label{normal}
\ee  
In order to state our result more efficiently we restrict attention to the surface $S:=S_+=(M,J_+)$
which may or may not be biholomorphic to $S_-$.
Recall also that the tautological twistor section $J_+$ identifies $S$ and $\Sigma_+ \subset Z$.
Similarly, the cycle $C_+ \subset \Sigma_+$ will be identified with $C \subset S$ while $\theta$
will denote the Lee form $\theta_+$.

\begin{Proposition} \label{leeb}
Let $(M,c,J_\pm)$ be an anti-self-dual bi-Hermitian compact surface with $b_1(M)$ odd.
Then the Lee form $\theta$ of the common Gauduchon metric on $S=(M,J_+)$ is harmonic
and using the above notations, the Lee bundle $\mathcal L _\theta$
is completely determined by the following equation:
\be \mathcal L ^2 _{\theta _{|_{C}}} = \nu^*_{C / \Sigma_{+}} \otimes  \nu_{C / \Sigma_{-}}
    \label{restrictedlee}   \ee
\end{Proposition}
\it Proof. \rm Under our assumptions  $b_1(M)=1$ so that the inclusion $i: C \hookrightarrow M$
induces an isomorphism of flat bundles 
$i_*:\operatorname{Hom}(\pi_1(M),\c))\to \operatorname{Hom}(\pi_1(C),\c))$
when $C$ is a cycle of rational curves
and an injection when it is a smooth elliptic curve.
In any case the Lee bundle $\mathcal L_\theta$ is completely determined by its restriction to $C$.

Now, by the main result of \cite[Theorem 2.1, Corollary 2.3]{po92} the Lee bundle of any
anti-self-dual Hermitian metric can be read off from the twistor space and setting
$\mathcal L = t^*(\mathcal L _{\theta})$ we have that the following equality of line bundles holds on $Z$:
$$X_+ = K^{-\frac{1}{2}}_Z  \otimes\mathcal L$$
similarly, because the Lee forms are opposite \cite[Lemma 3.4]{po97}
$$X_- = K^{-\frac{1}{2}}_Z   \otimes\mathcal L ^*$$ so that $$ \mathcal L ^2 = X_-^* \otimes X_+ $$ 
 restricting this equation to $C=\Sigma_+\cap\Sigma_-$ gives $\mathcal L ^2_{|_{C}} = (\Sigma^*_- \otimes \Sigma_+)_{|_{C}}$
 and therefore $\mathcal L ^2_{|_{C}}=\nu^*_+\otimes\nu_-$ by (\ref{normal}). Notice that each $\mathcal L^2 \in H^1(M,\r^{>0})$
 determines a unique $\mathcal L \in H^1(M,\r^{>0})$,  see (\ref{realflat}).          \hfill $\Box$  \bigskip
 
\begin{Remark} \rm As already mentioned, the above result applies to only three cases
\cite[lemma 3.2]{fp10}: Hopf surfaces, hyperbolic and parabolic Inoue surfaces.
When $S_+=(M,c,J_+)$ is an anti-self-dual Hopf surface it is a quotient of $\c^2 \setminus \{0\}$
by a diagonal action generated by a contraction of the form $(z,w) \mapsto (az,bw)$.
$S_\pm$ always have two elliptic curves,
denoted by $C$ and $E$, of periods $a$ and $b$ respectively.
In order to compute the Lee bundle $\mathcal L _\theta$ of the \lck metric on $S_+$
we notice that of course $ \nu_{C / \Sigma_{+}}$ induces the flat line bundle $C$ on $S_+$
which we associate to the non-zero complex number $a\in\c^* \cong \operatorname{Pic}^0(S_+)$.

The Hopf surface $S_-$ turns out to be the quotient by the conjugate contraction
$(z,w) \mapsto (\bar a z,\bar b w)$ \cite{po97} so that $ \nu_{C / \Sigma_{-}}$
corresponds to $\bar a \in \c^*$ and by Proposition \ref{minus} we get
that $\mathcal L^2 _\theta$ corresponds to the real line bundle $-|a|^2$.
We can exchange the roles of $C$ and $E$ in this case, leading to the reality condition
$|b|^2=|a|^2$ of \cite{po97}.
\end{Remark}

\begin{Remark}  \rm
The second case is the one of hyperbolic Inoue surfaces. In this case the complex structure
of $S_-$ turns out to be the \rm transpose surface  \it $^tS$ \cite{za01} \cite[\S 3.4]{fp10} of $S_+$.
For any of the two cycles $C$ of $S=S_+$ our Proposition above states that the Lee bundle
is induced by the tensor product of the conormal bundle of $C$ in $S$ with the normal
bundle of $C$ in $^tS$.
Now, the transpose hyperbolic Inoue surface $^tS$ is uniquely determined from $S$
by changing the cyclic order of the irreducible components of its maximal curves;
we can see in this way that the inclusion of either cycle $C$ induces opposite orientations
on $\pi_1(S)$ and this gives a geometric realization of Proposition \ref{minus}.
It also follows that the tensor product of the two (co)normal bundles
is always real, flat (as the two contributions coming from the first Chern class $c_1(C)$ cancel out)
and independent of the choice of $C$.
\end{Remark}

The most interesting case turns out to be the last case of parabolic Inoue surfaces.
In what follows we will identify the Lee class of \lck bi-Hermitian structures and clarify 
the reality condition appearing in the constructions of \cite{le91} and \cite{fp10}.

\begin{Proposition} \label{real}
Let $(M,c,J_\pm)$ be any anti-self-dual bi-Hermitian structure on a parabolic Inoue surface $S$.
Then the Lee form $\theta$ of the common Gauduchon metric is harmonic
and the Lee bundle satisfies $\mathcal L^*_\theta = C$
where $C$ is the unique cycle of rational curves.
It follows that  $C$ must represent a flat holomorphic line bundle  of real type.
A parabolic Inoue surface $S$ satisfying this condition will be called of ``real" type.
\end{Proposition}

\it Proof. \rm The irreducible components of the cycle $C$ of any Enoki surface,
in particular of a parabolic Inoue surface $S$, all have self-intersection equal to $-2$
and therefore its self-intersection $C^2=0$ showing that $C$ represents a flat line bundle always.
As already mentioned the inclusion induces an isomorphism $Pic^0(C)\cong Pic^0(S)$
sending the normal bundle $\nu_{C / \Sigma_{+}}$ to $C \in Pic^0(S)$.

 The notion of transposition does not apply in this case -- as all self-intersections are equal --
so that  $\Sigma_+$ and $\Sigma_-$ turn out to be isomorphic Inoue surfaces \cite{fp10}.
Since the two Lee forms are opposite, we get by Proposition \ref{minus} that
$\nu^*_{C / \Sigma_{-}}=C\in Pic^0(S)$ and therefore $(\mathcal L^*_\theta)^2 = C^2$.
It follows that $C=\mathcal L^*_\theta \in Pic^0_\r(S)$
therefore $C$ must be real - i.e. is a singular elliptic curve of real period.
\hfill $\Box$  \bigskip

\begin{Remark} \rm In \cite{fp10} we constructed families of bi-Hermitian structures on parabolic and hyperbolic Inoue surfaces with fixed complex structure, and varying Riemannian metric $g_t$. The real dimension of these families can be computed from the twistor space and turned out to be $b$. The above result shows that these  families of \lck metrics all have the same fixed Lee class which only depends on the complex structure $J$ and not on the metric $g_t$.
\end{Remark}

\section{New bi-Hermitian metrics with $\mathbf t =1$}

We start this section by recalling a strong result of Apostolov-Bailey-Dloussky

\begin{Theorem} \cite{abd15}  \label{abd}
Let $S=(M,J)$ be a compact surface in class $VII$ which admits a l.c.K. metric and
denote by $\mathcal L _{\theta}$  the Lee bundle of this metric. 
Suppose $S$ also admits a NAC divisor $\T=F\otimes K^*$ satisfying $F=\mathcal L _{\theta}$
then, $M$ has bi-Hermitian metrics $(g,J_\pm)$ such that the following holds:
\begin{enumerate}
\item $J_+ = J$ 
\item $\T=F\otimes K^*$ is the fundamental divisor of the bi-Hermitian structure $(M,[g],J_\pm)$ 
\item $\T$ is connected: $b_0(\T)=: \mathbf t =1$
\end{enumerate}
\end{Theorem}

\begin{Remark} \rm Assume as usual that $S$ is minimal,
the existence of the NAC divisor $\T$ forces $S$ to be a Hopf surface when $b_2(S)=0$ 
or else a Kato surface, for $b_2(S)>0$.
These surfaces always admit l.c.K. metrics \cite{go99} \cite{be00} \cite{br11}. 
So the important point here is the relation between NAC divisors (which are often unique on these surfaces) and Lee bundles.

Every Hopf surface turns out to admit \lck metrics with potential \cite{go99} and is a small deformation
of a diagonal Hopf surface admiting l.c.K. metrics with parallel Lee form, this was used in 
\cite[5.1]{ad16} to show that they have Lee bundles in every possible de Rham class.

For Kato surfaces the situation is different: 
because the universal cover has compact curves 
they cannot have \lck metrics with potential 
and the result of Proposition \ref{real} points out a new Lee class in parabolic Inoue surfaces 
of real type which, thanks to Theorem \ref{abd}, 
will be shown to produce bi-Hermitian metrics of a new type.
 \hfill $\bigtriangleup$
\end{Remark}

\medskip
Recall that all Kato surfaces have exactly $b_2>0$ rational curves some of which (maybe all) 
forming at least one (and at most two) cycle(s).
We say that a Kato surface $S$ is \it unbranched \rm if every rational curve belongs to a cycle.
The other Kato surfaces are \it branched \rm and they are also called \it intermediate. \rm

\medskip
Together with the results of \cite{fp10} the following gives a definite pictures of bi-Hermitian
Kato surfaces with no branches.

\begin{Theorem}  \label{parinoue}
An unbranched Kato surface $S$ admits bi-Hermitian metrics with $\mathbf t =1$ 
if and only if $S$ is a parabolic Inoue surface of real type.
\end{Theorem}
{\it Proof.} It is shown in  \cite{fp14} and \cite{abd15} that if $(M,c,J_\pm)$ is a class-VII$_0$
bi-Hermitian surface with fundamental NAC divisor $\T=F\otimes K^*$ and $\mathbf t =1$ then 
$S=(M,J_+)$ is either a Hopf  surface (which is not a Kato surface) 
or a parabolic Inoue surface in which case $K^*=E\otimes C$ with the following properties:
$E$ is the unique (smooth) elliptic curve, $C$ the unique cycle and they are disjoint.
This implies that the only possibility for $\T$ to be non-empty and connected is therefore $F=C^*$.

In particular, the cycle $C$ must correspond to a flat holomorphic line bundle of real type, 
i.e. the parabolic Inoue surface is of real type.

Viceversa, let $S$ be such a parabolic Inoue surface and let $(c,J_\pm)$ be
an anti-self-dual bi-Hermitian structure with $S=(M,J_+)$ as in the previous section --
in particular \lck
We showed in Proposition \ref{real} that the Lee bundle of such a metric satisfies 
$\mathcal L^*_{\theta} = C$, the conclusion then follows from Theorem \ref{abd}.
 \hfill$\Box$\medskip

\begin{Remark}  \rm
These are the first examples of bi-Hermitian Kato surfaces with $\mathbf t =1$.
In order to complete the classification it would now be enough to consider branched Kato surfaces
(also called \it intermediate \rm) of index $1$ 
-- i.e. admitting a divisor $D=G \otimes K^*$ with $c_1(G)=0$.
It is known in this case \cite{dl06} that the NAC divisor is unique and contains 
the maximal curve which is connected, in particular we always have $\mathbf t =1$
if a branched Kato surface is bi-Hermitian.
In \cite{fp15} we described some geometric obstructions to existence of bi-Hermitian
metrics on intermediate -- i.e. branched -- Kato surfaces; 
however,  it turns out that the techniques developed in this work are not enough to produce existence results for bi-Hermitian intermediate Kato surfaces. 
 \hfill $\bigtriangleup$      \end{Remark}

\section{Deformations}

The deformation theory of \lck structures with Lee form $\theta$ on a compact complex manifold
$X=(M,J)$ is governed by the twisted de Rham complex defined by the differential operator
$d_{\theta}:=(d-\theta\wedge)$ which satisfies $d_\theta^2=0$ because $\theta$ is closed.
We will denote by $H^r(M,L_{\theta}^*)$ the corresponding cohomology groups
which only depend on the smooth structure of $M$;
recall from the introduction that the fundamental $(1,1)$-form $\omega\in H^2(M,L_{\theta}^*)$.
It is also inportant to consider the twisted Dolbeaut complex defined by the differential operator
$\overline{\partial}_{\theta}:=(\overline{\partial}-\theta^{0,1} \wedge)$
whose cohomolgy groups $H^{p,q}(X,\mathcal L _{\theta}^*)$
are isomorphic to the sheaf cohomology groups
$H^q (X,\Omega^p \otimes \mathcal L_{\theta}^*)$ by Dolbeaut theorem
and therefore heavily depend on the complex structure $J$.

More precisely it is shown by Goto \cite{go14} that the deformation theory of \lck structures
is completely unobstructed  if the $\partial\overline{\partial}$-lemma holds at
degree $(1,2)$ and if furthermore  $H^3(M,L_{\theta}^*)=0$.
Notice that the first condition holds automatically when
$H^{0,2}(X,\mathcal L^*_\theta) \cong H^2(X,\mathcal L _{\theta}^*) = 0$
while the second is a smooth invariant one.

\begin{Remark} \rm A simple consequence of the above discussion and the following one,
is that if $\theta$ is an isolated Lee class on $S\in$ VII$_0$, then $S$ is either an unknown
surface without a cycle of rational curves or else it must be an Inoue-Bombieri surface (with  $b_2(S)=0$)
and furthermore the Lee class must represent the anticanonical bundle, see \cite{go14} and \cite{ot16}.     \hfill $\bigtriangleup$
\end{Remark}
\medskip

\begin{Remark} \rm Since the appearance of the first version of this paper \cite{fp16} several of the results in this section have been extended in \cite{ad16a} to surfaces with $b_1=1$.           \hfill $\bigtriangleup$
\end{Remark}
\bigskip

We are interested in \lck and bi-Hermitian structures on small deformations
of a given Kato surface $S$, for this reason let us notice that in dimension two we have
$H^2(S,\mathcal L _{\theta}^*) \cong H^0(S,K\otimes \mathcal L _{\theta})$ by Serre duality;
for this reason we now recall some known results about twisted pluri (anti-)canonical sections on
class-VII surfaces, see \cite{en81}, \cite{do99} and \cite{dl06}:

\begin{Proposition}\label{kato} Let $S$ be a surface in class $VII_0^+$ and
  $\mathcal L \in \operatorname{Pic}^0(S)$ then,
\begin{enumerate}
\item  $H^0(S,K^m \otimes \mathcal L) = 0$ for all $m>0$ and all
  $\mathcal L \in \operatorname{Pic^0}(S)$.
\item $H^0(S,\mathcal L) \neq 0$ for some $\mathcal L \neq \mathcal O$
 if and only if $S$ is an Enoki surface and $\mathcal L =C^m$
 is a positive power of the cycle of rational curves $C\subset S$.
\item  Assume $S$ is not an Enoki surface. Then,
         $H^0(S,K^{-m} \otimes \mathcal L) \neq 0$ for some $m>0$ if and only if $S$ is a Kato surface.
         The least such $m$ is called the index of the Kato surface and the line bundle
         $\mathcal L \in \operatorname{Pic}^0(S)\cong \c^*$ is unique.
\item When $S$ is an Enoki surface, $H^0(S,K^{-m} \otimes \mathcal L) \neq 0$ with $m>0$
          holds if and only if $S$ is parabolic Inoue in which case the index $m=1$. Furthermore
          $\mathcal L$ is either trivial or $\mathcal L=C^{-m}$ where $C$ is the unique cycle.
\end{enumerate}
\end{Proposition}

The following computation follows easily from vanishing results which have already been observed in
\cite{llmp03}, \cite[prop.6.1]{go14} and \cite{ad16}.

\begin{Lemma} \label{blowup}
Let $S$ be a Kato surface (or more in general a class-VII surface with a cycle of rational curves)
and $L$ a non-trivial flat line bundle on $S$, then $H^p(S,L) = 0$ for all $p\neq 2$ while dim$_\r H^2(S,L)=b_2(S)$.
\end{Lemma}
{\it Proof.}
The statement is invariant under diffeomorphisms therefore is enough to recall a result of Nakamura
\cite{na89} stating that $S$ is diffeomorphic to a blown up Hopf surface,  we denote  by $\check{S}$
the minimal model.
After a deformation, we can assume that $\check{S}$ is bi-Hermitian with $\t=2$ and admits
Vaisman metrics with arbitrary Lee form for both complex structures.
It then follows by a result of  \cite{llmp03} that $H^p(\check{S},L) = 0$
for every non-trivial flat line $L$.
Finally, the result follows because the blow up map $b:S\to \check{S}$ induces a surjection
$H^p( \check{S},L) \to H^p(S,b^*(L))$ for every $p\neq 2$, see for example
\cite[Lemma 2.3]{fp14} or  \cite[Lemma 3.3]{ad16}.
To conclude is enough to recall that by Hirzebruch-Riemann-Roch $\chi(S,L)=\chi(S)=b_2(S)$.
 \hfill$\Box$\medskip


 \smallskip

 \begin{Remark} \rm
 Notice that instead when $L$ is trivial, i.e. the differential operator is the usual de Rham operator, 
 we have Betti numbers $b_i(S)=1$ for each $i\neq 2$.        \hfill $\bigtriangleup$
 \end{Remark}
\medskip

 From the work of Goto \cite{go14} we then get the following result which was
also proved by Brunella in his explicit construction \cite[Theorem 2]{br10}.

\begin{Corollary}\label{unobstruction}
 The deformation theory of \lck structures on Kato surfaces is completely unobstructed in the
  sense that for a \lck Kato surface $S$ with Lee bundle $L$,
  every small deformation $(S_t,L_s)$ of the pair is still \lck with Lee bundle $L_s$.
\end{Corollary}

For a compact complex surface $S$ with odd first Betti number Apostolv-Dloussky
considered in \cite{ad16} the following set
$$\mathcal C (S) := \{ \theta\in H^1_{dR}(S,\r) \;|\; \theta
   \text{ is the Lee form   of a \lck metric} \}. $$
Of course when $b_1(S)=1$ we have $\mathcal C (S) \subset (-\infty,0) \subset H^1_{dR}(S,\r)$ because Lee bundles have negative degree, see \S 1 above.
When $S$ is a Kato surface $\mathcal C (S)$ is open and contains an interval $(-\infty, B_S)$ 
by Brunella's construction \cite[Remark 3]{br11};  the non-positive constant $B_S$ varies smoothly under deformations and
recall $B_S = 0$ when $S$ is a (blown up) Hopf surface \cite[Prop. 5.1]{ad16}.

\begin{Remark} \rm Brunella's construction shows that for any Kato surface -- i.e. for a fixed
complex structure $S=(M,J)$ on the smooth four-manifold $M=(S^1\times S^3) \# m \overline{\mathbb{CP}}_2$ --
there exist \lck metrics whose Lee form $B_S$ is close to  $-\infty \in H^1_{dR}(M,\r)$.
The result of Proposition \ref{real} can be used to show that for any given $\gamma\in (-\infty,0) \subset H^1_{dR}(M)$
there is a parabolic Inoue surface -- i.e. a particular complex structure $J$ on $M$ --
such that $\gamma$ is the Lee form of a \lck metric $g$ on $(M,J)$; the correspondence is given by taking the opposite of the $g$-volume of the cycle of rational curves which in turns is the trace tr$(S)$ of the parabolic Inoue surface and varies between $0$ and $1$, see \cite[Lemma 4.21]{dl14}.

For this reason, when $\gamma$ is arbitrarily  close to $0\in H^1_{dR}(M)$ these metrics seem to be far from the ones constructed by Brunella. The same holds by Corollary \ref{unobstruction} for any small deformation $(M,J_t)$ of parabolic Inoue surfaces. Because every Enoki surface is a small deformation of a parabolic Inoue surface and Enoki surfaces fill up a big open set in the moduli space of Kato surfaces, we see that $\gamma$ can be taken to be arbitrarily close to $0\in H^1_{dR}(M)$ 
on an open set in the moduli space of Kato surfaces.
 \hfill $\bigtriangleup$                        \end{Remark}

The above observations coupled with results of Apostolov-Dlousskly give the following:
\begin{Corollary} For any Enoki surface $(S,J)$ which is a small deformation of a ``real" parabolic Inoue surface, the set
$$\mathcal T (S):=\{\alpha\in H^1_{dR}(S,\r) \; | \; \alpha \text{ is the Lee class of a LCS structure taming } J \}$$
of Lee classes of locally conformally symplectic structures (LCS) taming the complex structure is an open set $\mathcal T (S)\subset (-\infty,0) \subset H^1_{dR}(S,\r)$ containing the closed interval $ (-\infty,\gamma]$.
Here $\gamma$ is a negative multilple of $\log | trace(S) |$, it goes to $0$ as $| trace(S) |$ tends to $1$.
\end{Corollary}
{\it Proof.} Start with a parabolic Inoue surface $S$ and apply  \cite[Proposition 4.4]{ad16} with $b\in\mathcal C (S)$ a \lck Lee form of Brunella (arbitrarily close to $-\infty$) and $c\in\mathcal G (S)$ the harmonic Lee form of our bi-Hermitian metrics to conclude that $(-\infty,c]\subset\mathcal T (S)$ because $H^2(S,\mathcal L _c)=0$ by (\ref{kato}) and Serre duality. The same result holds for Enoki surfaces since they are small deformations of parabolic Inoue ones and since the stability result of (\ref{unobstruction}) holds in the more general context of LCS structures yielding in particular that $\mathcal T (S)$ is always open for Kato surfaces, \cite[Theorem 1.3]{ad16a}.
 \hfill$\Box$    \medskip

{\bf Concluding remarks.} \it Intermediate surfaces. \rm 
Large families of deformations of Kato surfaces with second Betti number $b$ have been introduced and studied by Dloussky in \cite{dl14} and \cite{dl16}. Over the $b$-dimensional base space $\mathcal B$ of these deformations, Enoki surfaces fill up an open subset and the boundary $\partial\mathcal B$ is a union of strata of different dimensions ranging from $(b-1)$ down to $0$. The $0$-dimensional stratum consists of hyperbolic and half Inoue surfaces which are rigid under deformations. The positive dimensional strata of $\partial\mathcal B$ consist instead of so called intermediate surfaces. 
If we denote by $\gamma\in (-\infty, 0)\subset H^1_{dR}$ the Lee class of an anti-self-dual bi-Hermitian metric on real parabolic Inoue surfaces, then $\gamma$ is still a Lee class on nearby deformations which are necessarily Enoki surfaces; from \cite[Lemma 4.21]{dl14} it follows that $\gamma$ goes to $-\infty$ when approaching the boundary $\partial\mathcal B$ -- because the trace goes to 0 -- and therefore our results fail to produce new Lee classes (i.e. different from Brunella's construction) on intermediate surfaces.

On the other hand, let us consider bi-Hermitian structures on hyperbolic Inoue surfaces; these surfaces belong to the $0$-dimensional stratum in a large family of Kato surfaces and we would like to investigate existence of bi-Hermitian structures on nearby intermediate surfaces; every such surface (of index $1$) has a unique NAC divisor $F \otimes K^*$ and is therefore a good candidate to admit bi-Hermitian metrics when the Gauduchon degree of $F$ satisfies deg$F<0$. The flat line bundle $F$ is controlled by the  leading coefficient of the Favre germ of the surface, however  it follows from \cite[4.24]{dl14} that deg$F$ tends to $+\infty$ when the surface approaches a lower dimensional stratum and the conclusion is that there is no bi-Hermitian structure on intermediate Kato surfaces which are small deformations of  hyperbolic Inoue ones.
 \hfill $\bigtriangleup$

\newcommand{\bysame}{\leavevmode\hbox to3em{\hrulefill}\,}


\begin{thebibliography}{99}
\bibitem[Ap01]{ap01} Apostolov V.,
    Bihermitian surfaces with odd first Betti number.
    {\em Math. Zeit.} {\bf 238} (2001), 555-568.
\bibitem[ABD15]{abd15}    Apostolov V., Bailey M., Dloussky G., 
  From locally conformally K\"ahler to bi-Hermitian structures on non-K\"ahler complex surfaces. 
  {\it Math. Res. Lett.} {\bf 22} (2015) 317--336.
\bibitem[AD08]{ad08} Apostolov V., Dloussky G.,
    Bihermitian metrics on Hopf surfaces.
    {\it Math. Res. Lett.} {\bf 15} (2008), no. 5, 827--839.
\bibitem[AD16]{ad16} \bysame 
  Locally conformally symplectic structures on compact non-K\"ahler complex surfaces.
 {\it Int. Math. Res. Not. IMRN} (2016), no. 9. 2717--2747.
 \bibitem[AD16a]{ad16a} \bysame
 On the Lee classes of locally conformally symplectic complex surfaces.
 arXiv:1611.00074  
\bibitem[AGG99]{agg99} Apostolov, V., Gauduchon, P., Grantcharov, G.,
    Bi-Hermitian structures on complex surfaces.
    {\em Proc. London Math. Soc.} (3) {\bf 79} (1999), no. 2, 414--428;
    {\em Proc. London Math. Soc.} (3) {\bf 92} (2006), no. 1, 200--202.
\bibitem[AG07]{ag07} Apostolov V., Gualtieri M.,
    Generalized K\"ahler manifolds, commuting complex structures
    and split tangent bundle.
    {\em Comm. Math. Phys.} {\bf 271} (2007) 561-575.
\bibitem[AHS78]{ahs} Atiyah M.F., Hitchin N.J., Singer I..M.,
    Self-duality in four-dimensional Riemannian geometry.
    {\em Proc. R. Soc.} A {\bf 362}  (1978) 425-461.  
\bibitem[Be00]{be00} Belgun, F. A.,
    On the metric structure of non-K\"ahler complex surfaces.
    {\em Math. Ann.} {\bf 317} (2000) 1--40.
\bibitem[Bo88]{bo88} Boyer, C.P.,
    A note on hyperhemitian four-manifolds.
    {\em Proc. Am. Math. Soc.} {\bf 102} (1988) 157-164.
\bibitem[Br10]{br10} Brunella, M.,
 Locally conformally K\"ahler metrics on certain non-K\"ahlerian surfaces.
{\it Math. Ann.} (2010) {\bf 346} 629 -- 639.
\bibitem[Br11]{br11} \bysame
   Locally conformally K\"ahler metrics on Kato surfaces.
   {\em Nagoya Math. J.} {\bf 202} (2011), 77-81.
\bibitem[CG11]{cg11} Cavalcanti G., Gualtieri M.,
  Blowing up generalized K\"ahler 4-manifolds. 
  {\it Bull. Braz. Math. Soc. (N.S.)} {\bf 42} (2011), no. 4, 537--557. 
  \bibitem[Dl06]{dl06} Dloussky, G.,
    On surfaces of class  VII$_0$ with numerically anticanonical divisor. 
    {\em Amer. J. Math.} {\bf 128} (2006), no. 3, 639--670.
\bibitem[Dl14]{dl14} \bysame
  From non-K\"ahlerian surfaces to Cremona group of $\cp_2$.  
    {\it Complex manifolds} {\bf 1} (2014) 1 -- 33. 
\bibitem[Dl16]{dl16} \bysame
  Special birational structures on non-K\"ahlerian complex surfaces.
  {\it J. Math. Pures Appl. (9) } {\bf 106} (2016) 76--122. 
\bibitem[DO99]{do99} Dloussky, G. ; Oeljeklaus, K.
  Vector fields and foliations an compact surfaces of class VII$_0$.
  {\it Ann. Inst. Fourier} {\bf 49} (1999) 1503 -- 1545. 
\bibitem[DOT03]{dot03} Dloussky, G.; Oeljeklaus, K.; Toma, M. 
 Class VII$_0$ surfaces with $b_2$ curves. 
 {\it Tohoku Math. J.} (2) {\bf 55} (2003) 283-Ð309.
\bibitem[En81]{en81}  Enoki, I. 
   Surfaces of class VII$_0$ with curves. 
   {\em T\^ohoku Math. J.} (2) {\bf 33} (1981)  453Ð492.
\bibitem[FP10]{fp10} Fujiki, A., Pontecorvo, M.,
    Anti-self-dual bihermitian structures on Inoue surfaces.
    {\em J. Differential Geom.} {\bf 85} (2010), no. 1, 15--71.
\bibitem[FP14]{fp14} \bysame  
   Twistors and bi-Hermitian surfaces of non-K\"ahler type.
   {\it Symmetry, Integrability and Geometry: Methods and Applications (SIGMA).}
   {\bf 10} (2014) Paper 042, 13 pp.
\bibitem[FP15]{fp15} \bysame 
  Numerically anti-canonical divisors on Kato surfaces.
  {\it J.  Geom. Phys.} {\bf 91} (2015) 117--130.    
\bibitem[FP16]{fp16} \bysame
Bi-Hermitian metrics on Kato surfaces.  {\tt arXiv:1607.00192v1}.
\bibitem[Ga81]{ga81} Gauduchon P.,
  Le th\'eor\`eme de dualit\'e pluriharminique.
  {\it C.R. Acad. Sci Paris} {\bf 293} (1981) 59--63.
\bibitem[Ga84]{ga84} \bysame  
   La $1$-form de torsion d'une vari\'et\'e hermitienne compacte.
    {\em Math. Ann}  {\bf 267} (1984) 495--518.
\bibitem[GO99]{go99} Gauduchon P., Ornea L.,
    Locally conformally K\"ahler metrics on Hopf surfaces.
    {\em Ann. Inst. Fourier (Grenoble)} {\bf 48} (1998) 1107-1127.    
\bibitem[Go12]{go12} Goto, R.,
  Unobstructed K-deformations of generalized complex structures and bi-Hermitian structures. 
  {\it Adv. Math.} {\bf 231} (2012) 1041--1067.
\bibitem[Go14]{go14} \bysame  
  On the stability of locally conformal K\"ahler structures. 
  {\it J. Math. Soc. Japan} {\bf 66} (2014) 1375--1401. 
\bibitem[Hi06]{hi06}  Hitchin, N. J.,
    Instantons, Poisson structures and generalized K\"ahler geometry.
    {\em Commun. Math. Phys.} {\bf 265} (2006) 131--164.
\bibitem[Hi07]{hi07} \bysame
    Bihermitian Metrics on Del Pezzo surfaces.
    {\it J. Symplectic Geom.} {\bf 5} (2007) 1--8.        
\bibitem[Ka77]{ka77} Kato, Ma.,
   Compact complex manifolds containing ``global'' spherical shells.
    {\it I. Proceedings of the International Symposium on Algebraic Geometry}
    (Kyoto Univ., Kyoto, 1977), pp. 45--84, Kinokuniya Book Store, Tokyo, 1978.
\bibitem[Le91]{le91} Lebrun, C.R.,
    Anti-self-dual Hermitian metrics on blown-up Hopf surfaces.
    {\em Math. Ann.} {\bf 289} (1991) 383--392.
\bibitem[dLLMP03]{llmp03} M. de Le\'on, B. L\'opez, J.C. Marrero, E. Padr\'on.
On the computation of the Lichnerowicz-Jacobi cohomology.
{\it J. Geom. Phys.} {\bf 44} (2003) 507 -- 522.
\bibitem[Na89]{na89}  Nakamura, I. 
  Classification of non-KŠhler complex surfaces. (Japanese) Translated in 
  {\it Sugaku Expositions} {\bf 2} (1989) 209--229. 
  Sugaku 36 (1984), no. 2, 110Ð124.
\bibitem[Ot16]{ot16} Otiman A.,
      Morse-Novikov cohomology of locally conformally K\"ahler surfaces.
     {\it  Math. Z.} (2018). https://doi.org/10.1007/s00209-017-1968-y. 
\bibitem[Po92]{po92} Pontecorvo M.,
    On twistor spaces of anti-self-dual Hermitian surfaces.
    {\it Trans. Amer. Math. Soc.} {\bf 331} (1992) 653--661.
\bibitem[Po97]{po97} \bysame  
    Complex structures on Riemannian four-manifolds.
    {\em Math. Ann.} {\bf 309} (1997), no. 1, 159--177.
\bibitem[Sa09]{sa09} Salamon, S.,
Complex structures and conformal geometry. 
{\it Boll. Unione Mat. Ital.} (9) {\bf 2} (2009) 199 -- 224. 
\bibitem[Tr74]{tr74} Tricerri, F.,
 Sulle variet\`a dotate di due strutture quasi complesse
    linearmente indipendenti.
    {\em Riv. Mat. Univ. Parma} (3) {\bf 3} (1974) 349-358.
\bibitem[Za01]{za01} Zaffran D., Serre problem and Inoue-Hirzebruch surfaces. 
  {\it Math. Ann.} {\bf 319} (2001), no. 2, 395--420. 
\end{thebibliography}
\end{document}